\documentclass[10pt,a4paper]{amsart}
\usepackage{amssymb}
\usepackage{amsmath}
\usepackage{hyperref}
\usepackage{amsfonts,a4wide}
\usepackage{bbm}
\usepackage{enumerate}
\usepackage{inputenc}
\usepackage{shadow}
\usepackage[all]{xy}
\usepackage{tikz-cd}
\usepackage{cleveref}
\usepackage{mathdots}
\usepackage[fleqn]{mathtools}
\usepackage[a4paper, bindingoffset=0.2in, left=2.15cm, right=2.15cm, top=3.5cm, bottom=3.5cm, footskip=.25in]{geometry}

\theoremstyle{definition}
\newtheorem{thm}[subsection]{Theorem}
\newtheorem{defn}[subsection]{Definition}
\newtheorem{prop}[subsection]{Proposition}

\newtheorem{cor}[subsection]{Corollary}
\newtheorem{lemma}[subsection]{Lemma}

\newtheorem{remark}[subsection]{Remark}
\newtheorem{example}[subsection]{Example}

\newcommand{\R}{\mathbb{R}}

\newcommand{\CP}{\mathbb{CP}}

\newcommand{\GL}{\operatorname{GL}}
\newcommand{\Z}{\mathbb{Z}}

\newcommand{\JJ}{\mathcal J}
\newcommand{\C}{\mathbb{C}}

\newcommand{\HH}{\mathbb{H}}

\newcommand{\la}{\mathfrak a}
\newcommand{\im}{\operatorname{Im}}

\newcommand{\so}{\mathfrak{so}}

\newcommand{\uu}{\mathfrak u}
\newcommand{\lp}{\mathfrak p}

\DeclareMathOperator{\Sp}{Sp}

\DeclareMathOperator{\Hom}{Hom}

\DeclareMathOperator{\id}{id}
\DeclareMathOperator{\Gr}{Gr}

\DeclareMathOperator{\Id}{Id}

\DeclareMathOperator{\End}{End}

\DeclareMathOperator{\Span}{Span}
\DeclareMathOperator{\grph}{graph}

\begin{document}

\title{Intersections of complex structures}

\author{Gustavo Granja}
\address{Center for Mathematical Analysis, Geometry and Dynamical Systems, Instituto Superior T\'ecnico,
Universidade de Lisboa, Av. Rovisco Pais, 1049-001 Lisboa, Portugal}
\email{gustavo.granja@tecnico.ulisboa.pt}

\author{Aleksandar Milivojevi\'c}
\address{University of Waterloo,
Faculty of Mathematics,
200 University Avenue West
Waterloo, ON, Canada N2L 3G1}
\email{amilivojevic@uwaterloo.ca}


\begin{abstract}
    We study the sets of planes in an even 
    dimensional real vector space $V$ which are simultaneously stabilised by a pair of complex structures on $V$. We completely describe these sets of planes for pairs of orthogonal complex structures. Generically,
    the number of such planes is finite. We compute this number for orthogonal complex structures and prove that it gives a lower bound for the number of planes simultaneously stabilised by a generic pair of complex structures on $V$.
\end{abstract}

\maketitle

\setcounter{tocdepth}{1}
\tableofcontents

\section{Introduction}

In previous work on the topology of the space of almost complex structures on six--manifolds (see \cite{FGM}), an important role was played by a certain parametrisation, due to Battaglia \cite{Ba}, of the space $J(\R^6)$ of orthogonal linear\footnote{A linear complex structure on a real vector space $V$ is an endomorphism $J\colon V \to V$
which squares to $-\Id_V$. In this paper, we deal only with linear complex structures, so 
from now on we will omit the adjective ``linear''.} complex structures on $\R^6$ compatible with the standard orientation. The validity of this parametrisation was found in \cite[Proposition 2.1]{FGM} to be equivalent to the fact that any two distinct elements of $J(\R^6)$ agree on exactly one common complex line.
We describe this fact by saying that two distinct orthogonal complex
structures on $\R^6$ compatible with the orientation \emph{intersect} 
on a line. 

The aim of this note is to address the natural follow-up question of 
describing the sets of $2k$-planes on which two complex structures on 
$\R^{2n}$ intersect. For $k>1$ it is no longer reasonable to expect that two
complex structures will actually agree on a $2k$-plane, so the notion of
intersection at a plane $P$ must be relaxed to the condition that 
$P$ is a complex plane with respect to both complex structures. 
In the future we hope to use the results of this paper to 
study the topology of the space of almost complex structures on higher 
dimensional manifolds. 

Geometrically, the statement that two complex structures $J_0$ and $J_1$ intersect
at a $2k$-plane $P$ can be phrased in the following way: 
$J_0$ and $J_1$ determine embedded copies $\Gr^{J_i}_{2k}(\R^{2n})$
of the Grassmannian of complex $k$-planes in the real Grassmannian $\Gr_{2k}(\R^{2n})$, 
namely, the subsets consisting of planes stabilised 
by the $J_i$; a plane $P$ where $J_0$ and $J_1$ intersect is simply a point in the 
intersection $\Gr_{2k}^{J_0}(\R^{2n}) \cap \Gr_{2k}^{J_1}(\R^{2n})$.

The problem of counting intersection points can be approached homologically: the fundamental classes of the complex Grassmannians are middle-dimensional homology classes in the real Grassmannian and so
have a homological intersection number. 
We will consider the intersection in the Grassmannian $\Gr_{2k}^+(\R^{2n})$ of oriented $2k$-planes so as to be able to distinguish planes
on which $J_0$ and $J_1$ induce the same orientation from those 
on which they induce opposite orientations. 
Our first main result is a computation of these homological intersection numbers.
For $0\leq k \leq n$ let 
\begin{equation*}
\sigma(k,n)=\begin{cases}
    0 & \text{ if } k \text{ is odd and } n \text{ is even,}\\
    \dbinom{ \left \lfloor{\frac{n}{2}}\right \rfloor}{ \left \lfloor{\frac{k}{2}}\right \rfloor} & \text{ otherwise. }
\end{cases}
\end{equation*}
(see Table 1 below). If $k>n$ or $k<0$ we set $\sigma(k,n)=0$.
\begin{thm}
\label{main}
Let $J_0$ and $J_1$ be complex structures on $\R^{2n}$ and $1\leq k \leq n$.
\begin{enumerate}[(i)]
\item If $J_0$ and $J_1$ induce the same orientation on $\R^{2n}$, then  
$$[\Gr_{2k}^{J_0}(\R^{2n})] \cdot  [\Gr_{2k}^{J_1}(\R^{2n})] = \sigma(k,n) 
\in H_0(\Gr_{2k}^+(\R^{2n})).$$
\item If $J_0$ and $J_1$ induce opposite orientations on $\R^{2n}$, then 
$$[\Gr_{2k}^{J_0}(\R^{2n})] \cdot  [\Gr_{2k}^{J_1}(\R^{2n})] = \sigma(k,n-1)
\in H_0(\Gr_{2k}^+(\R^{2n})).$$
\end{enumerate}
\end{thm}

A more complete version of this theorem, which considers also planes of intersection where 
the complex structures induce opposite orientations, is proved below as Theorem \ref{main1} 
and Proposition \ref{oppor}.

Generically, the submanifolds $\Gr_{2k}^{J_0}(\R^{2n})$ and $\Gr_{2k}^{J_1}(\R^{2n})$
intersect transversally (see Proposition \ref{generic}), in which case the homological
intersection numbers described above are given by 
a signed count of the actual intersection points.
Hence Theorem \ref{main} gives a lower bound for the 
number of planes of intersection of a generic pair of complex structures (on which 
the complex structures induce the same orientation) - see Corollary \ref{estints}. 

When the complex structures are both orthogonal (with respect to some inner product)
we can completely describe the intersections (see Theorem \ref{descint}). Indeed, 
$\Gr_{2k}^{J_0}(\R^{2n}) \cap \Gr_{2k}^{J_1}(\R^{2n})$ depends only on the isomorphism 
type of the pair $(J_0,J_1)$ and, in the orthogonal case, these isomorphism classes 
correspond to the orbits of the action of $U(n)$ on $O(2n)/U(n)$. The basic theory
of symmetric spaces allows us to give a simple description of all possibilities
(see Proposition \ref{class}) and prove the following result (see Corollaries \ref{genorthopos} and \ref{genorthoneg}): 

\begin{thm}
\label{mainno}
Let $J_0,J_1$ be orthogonal complex structures on $\R^{2n}$. Generically, the number
of planes on which $J_0$ and $J_1$ intersect and induce the same orientation is the
homological intersection number described in Theorem \ref{main}.
\end{thm}

\subsection{Organization of the paper}
In section 2, we prove Theorem \ref{main} by making use of the Atiyah-Bott and 
Berline-Vergne equivariant localisation theorem. We also explain how these 
homological intersection numbers give a generic lower bound for the number of 
planes of intersection. 
In section 3, we classify pairs of orthogonal 
complex structures on finite dimensional Euclidean vector spaces. In section 4 we use 
this classification to completely describe the possible intersections of two 
orthogonal complex structures. Finally, in section 5, we discuss some ways in which
the intersections of general complex structures differ from intersections of orthogonal 
ones. 

\subsection{Acknowledgements} This paper originated in our attempt to compute 
certain Chern numbers which equal the intersection numbers in Theorem \ref{main}.
We were greatly helped in this task by Daniel Dugger and, especially, Nicholas Addington
who provided us with relevant Macaulay2 code that suggested the statement of Theorem \ref{main}.
We would also like to thank  Scott Wilson and Luis Fernandez for useful discussions, and 
Michael Albanese for detailed comments and 
suggestions that have improved the clarity and accuracy 
of the paper. Finally, we thank J. C. 
Alvarez-Paiva for his mathoverflow answer 
\url{https:// mathoverflow.net/questions/91583/the-rank-of-a-symmetric-space} 
which played an important role in the development of our paper.
The first author was partially supported by FCT/Portugal through 
CAMGSD, IST-ID, projects UIDB/04459/2020 and UIDP/04459/2020.

\section{Homological estimate for the intersection numbers}

Let $\Gr_{2k}(\R^{2n})$ denote the Grassmannian of $2k$-planes in $\R^{2n}$ and 
$\Gr^+_{2k}(\R^{2n})$ denote the Grassmannian of oriented $2k$-planes. Let 
$$
\Gr^+_{2k}(\R^{2n}) \xrightarrow{\pi} \Gr_{2k}(\R^{2n})
$$
denote the canonical double cover which forgets the orientation. We write 
$$\Gr^+_{2k}(\R^{2n}) \xrightarrow{\tau} \Gr^+_{2k}(\R^{2n})$$
for the involution which switches the orientation of each plane.

We consider the standard orientation for $\R^{2n}$. The standard charts
for $\Gr_{2k}(\R^{2n})$ (or $\Gr_{2k}^+(\R^{2n})$) at a plane $P$ are obtained by 
choosing $Q$ such that $P \oplus Q = \R^{2n}$ and then assigning to a plane $P'$
near $P$ the unique element of $\Hom(P,Q)$ which has $P'$ as its graph. 
Given ordered bases $(v_1,\ldots, v_{2k})$ for $P$ and 
$(w_1, \ldots, w_{2n-2k})$ for $Q$ we consider the ordered basis 
\begin{equation}
\label{basis}
(v_1^*\otimes w_1, \ldots, v_1^*\otimes w_{2n-2k}, v_2^* \otimes w_1, \ldots, v_{2k}^*
\otimes w_{2n-2k})
\end{equation}
for $\Hom(P,Q) = P^*\otimes Q$.

An orientation for $P$ together with the orientation of $\R^{2n}$ determine
an orientation for the complement $Q$. If $(v_1,\ldots, v_{2k})$ and 
$(w_1, \ldots, w_{2n-2k})$ are ordered 
bases for $P$ and $Q$ compatible with the orientations
we give
$$
T_P \Gr_{2k}^+(\R^{2n}) = \Hom(P,Q) = P^*\otimes Q
$$
the orientation determined by the basis \eqref{basis}. As $P$ and $Q$ are even
dimensional, the orientation of such a basis is actually independent of any choices: two
automorphisms $\phi \colon P \to P$ and $\psi \colon Q \to Q$ act on \eqref{basis} via
$(\phi^{-1})^*\otimes \psi$, and this automorphism of $P^*\otimes Q$ has determinant
$(\det \phi^{-1})^{2n-2k}(\det \psi)^{2k}>0$.

An element $g \in \GL(2n,\R)$ acts on $\Gr_{2k}^+(\R^{2n})$ by sending $P$ to the 
plane $g(P)$ with the orientation induced by the orientation on $P$. If we choose
$g(Q)$ as a complement for $g(P)$, $g$ induces the map $\Hom(P,Q) \to \Hom(g(P),g(Q))$ defined by the expression
$$
A \mapsto g A g^{-1}.
$$
If we pick bases for $P$, $Q$, and their images under $g$ for $g(P)$ and $g(Q)$, 
in the local coordinates determined by \eqref{basis} around $P$ and $g(P)$, the map
$g$ is the identity. In particular, the action of $\GL(2n,\R)$ preserves the orientation 
of $\Gr_{2k}^+(\R^{2n})$.

The involution $\tau$ which switches the orientation of a plane $P\in \Gr_{2k}^+(\R^{2n})$
induces the identity map on $\Hom(P,Q)$ and hence also preserves the orientation of
$\Gr_{2k}^+(\R^{2n})$. We conclude that the actions of $GL(2n,\R)$ and 
$\tau$ preserve the intersection product on the homology of $\Gr_{2k}^+(\R^{2n})$:

\begin{lemma}
\label{involu}
Let $\alpha,\beta$ be homology classes in $H_*(\Gr_{2k}^+(\R^{2n}))$
and let $g \in \GL(2n,\R)$. Then 
$$
\tau_*(\alpha) \cdot \tau_*(\beta) = \alpha \cdot \beta,
\quad \quad
g_*(\alpha) \cdot g_*(\beta) = \alpha \cdot \beta.
$$
\end{lemma}
\begin{proof}
An oriented map of manifolds preserves the intersection product on homology.
\end{proof}

Let 
$$\JJ(\R^{2n}) = \{ J \in \End(\R^{2n}) \colon J^2=-\Id\} \cong \GL(2n,\R)/\GL(n,\C)$$
denote the space of complex structures on $\R^{2n}$. An element $J\in \JJ(\R^{2n})$
determines for each $0\leq k \leq n$ a submanifold
$$
\Gr_{2k}^J(\R^{2n}) = \{ P \in \Gr_{2k}(\R^{2n}) \colon J(P)=P\} 
$$
of $\Gr_{2k}(\R^{2n})$. Assigning to a plane $P$ the orientation induced by $J$
we obtain a submanifold of $\Gr_{2k}^+(\R^{2n})$ which we denote in the same way.

These submanifolds are embedded copies of the complex Grassmannian $\Gr_k(\C^n)$ (and
are therefore canonically oriented) as they are the images under the projection 
onto the second coordinate of the fibers in the natural fiber bundle 
\begin{equation}
\label{bundle}
  E = \{ (J,P) \in \JJ(\R^{2n}) \times \Gr_{2k}^+(\R^{2n}) \colon J(P)=P,
J \text{ induces given orientation on } P\} \xrightarrow{\pi_1} \JJ(\R^{2n})
\end{equation}

(note that $E \cong \GL(2n,\R) \times_{\GL(n,\C)} \Gr_k(\C^n)$). A path in 
$\JJ(\R^{2n})$ connecting $J$ and $J'$ in $\JJ(\R^{2n})$ determines a homology 
between the cycles $\Gr_{2k}^J(\R^{2n})$ and $\Gr_{2k}^{J'}(\R^{2n})$ in 
$\Gr_{2k}^+(\R^{2n})$.

Finally, note that if $J(P)=P$, then we can pick a complement $Q$ for $P$ such 
that $J(Q)=Q$. In the local chart $\Hom(P,Q)$ for $\Gr_{2k}^+(\R^{2n})$ 
at $P$ the complex Grassmannian $\Gr_{2k}^J(\R^{2n})$ becomes the subspace
$$
\Hom_{\C}^J(P,Q) \subset \Hom(P,Q)
$$
of complex linear maps with respect to the complex structure determined by $J$
on $P$ and $Q$. The normal bundle at $P$ is represented by the space 
$\Hom_{\overline \C}^J(P,Q)$
of $J$-conjugate linear maps from $P$ to $Q$ (i.e. maps that anti-commute with the 
action of $J$ on $P$ and $Q$). 

Giving $\Hom(P,Q)$ the structure of a complex
vector space via post-composition with $J$, the decomposition 
\begin{equation}
\label{orinor}
\Hom(P,Q) = \Hom_\C^J(P,Q) \oplus \Hom_{\overline \C}^J(P,Q)
\end{equation}
is one of complex vector spaces. The orientation we have picked for $\Hom(P,Q)$
agrees with the orientation induced by the complex structure and hence is compatible 
with the canonical orientations for $\Hom_\C^J(P,Q)$ and $\Hom_{\overline \C}^J(P,Q)$.

Our aim in this section is to compute the intersection products
$$
[\Gr_{2k}^J(\R^{2n})] \cdot [\Gr_{2k}^{J'}(\R^{2n})], 
\quad [\Gr_{2k}^J(\R^{2n})] \cdot \tau_*[\Gr_{2k}^{J'}(\R^{2n})] 
\in H_0 (\Gr_{2k}^+(\R^{2n}))
$$
for arbitrary $J,J' \in \JJ(\R^{2n})$. These are 
homological counts of the number of planes $P$ stabilised by both $J$ and $J'$ 
on which $J$ and $J'$ induce the same, respectively the opposite, orientation.

As the homology class of $[\Gr_{2k}^J(\R^{2n})]$ depends only on the connected component
of $J$ in $\JJ(\R^{2n})$ and for any $g\in \GL(2n,\R)$ we have 
$$g_*[\Gr_{2k}^J(\R^{2n})]=[\Gr_{2k}^{gJg^{-1}}(\R^{2n})]$$
as well as $\tau g=g\tau$, using Lemma \ref{involu} we see that 
it is sufficient to perform the computations when 
$J=J'$ is the standard complex structure $J_0$ on $\R^{2n}$ (we then write 
$\Gr_k(\C^n) = \Gr_{2k}^{J_0}(\R^{2n})$) and when $J=J_0$ and $J'$ is any 
complex structure in the component of $\JJ(\R^{2n})$ not containing $J_0$.
We will consider these two cases in turn. 

Before this we observe that, generically, the submanifolds $\Gr_{2k}^J(\R^{2n}), 
\Gr_{2k}^{J'}(\R^{2n})$ intersect transversally at a finite number of points.
\begin{prop}
\label{generic}
The set of $(J,J') \in \JJ(\R^{2n})^2$ such that $\Gr_{2k}^J(\R^{2n})$ and 
$\Gr_{2k}^{J'}(\R^{2n})$ intersect transversally at a finite number of points
contains an open dense set. The analogous statement holds for $\Gr_{2k}^J(\R^{2n})$ and $\tau(\Gr_{2k}^{J'}(\R^{2n}))$.
\end{prop}
\begin{proof}
The projection onto the first coordinate on the set of $(J,J')$ satisfying the condition 
in the statement is a fiber bundle over $\JJ(\R^{2n})$, equivariant for the action induced 
by the diagonal action of $\GL(2n,\R)$ on $\JJ(\R^{2n})^2$ and the standard 
action on $\JJ(R^{2n})$. Thus we may assume that $J$ is the standard complex 
structure $J_0$. 

Let $g\in \GL(2n,\R)$ be such that $J'=gJg^{-1}$ so that 
$\Gr_{2k}^{J'}(\R^{2n})=g\Gr_{2k}^{J}(\R^{2n})$. The projection onto the second 
coordinate from the total space of the bundle \eqref{bundle}
can be identified with the map
$$
\GL(2n,\R)\times_{\GL(n;\C)} \Gr_k(\C^n) \xrightarrow{F} \Gr_{2k}^+(\R^{2n})
$$
determined by the formula $F(g,P)=gP$. This map is induced from a map 
$\GL(2n,\R)\times \Gr_k(\C^n) \to \Gr_{2k}^+(\R^{2n})$, given by the same expression,
which is clearly a submersion (its restriction to each subspace $\GL(2n,\R)\times \{P\}$
is a fiber bundle). It follows that $F$ is also a submersion.

By the parametric transversality theorem (applied to local trivializations of the bundle 
\eqref{bundle}), for an open dense set of $J' \in \JJ(\R^{2n})$ the inclusion
$$
\Gr_{2k}^{J'}(\R^{2n}) \hookrightarrow \Gr_{2k}^+(\R^{2n})
$$
(which is given by $P\mapsto F(g,g^{-1}P)$) 
is transverse to $\Gr_k(\C^n)$. As the manifolds $\Gr_{2k}^{J}(\R^{2n})$ have half the 
dimension of $\Gr_{2k}^+(\R^{2n})$ and  are compact, this means that they intersect 
transversally at a finite set of points.

The second statement follows from the argument above upon composing the map $F$ with 
the diffeomorphism $\tau$.
\end{proof}

\begin{remark}
\label{genericrem}
The previous statement admits several variations. For instance, 
we may replace $\GL(2n,\R)$ with the orthogonal group $O(2n)$ and arbitrary complex 
structures with orthogonal ones. 
\end{remark}

For transverse intersections, the homological intersection number can be computed as a 
signed count of the actual intersection points, therefore Proposition \ref{generic} 
implies that the homological intersection numbers we are about to compute provide
lower bounds for the actual number of intersection points between generic pairs of
complex structures.

Recall that we have defined for 
$0\leq k \leq n$ the quantity
$$
\sigma(k,n) = \begin{cases} 
0 & \text{ if } k \text{ is odd and } n \text{ is even,} \\ \dbinom{ \left \lfloor{\frac{n}{2}}\right \rfloor}{ \left \lfloor{\frac{k}{2}}\right \rfloor} & \text{ otherwise. }
\end{cases}
$$
Table 1 describes $\sigma(k,n)$ for low values of $k$ and $n$.
\begin{table}[htbp]\label{values}
\centering
\begin{tabular}{ |c||c | c |  c| c | c | c | c | c | c | c| c | c|c|c|c| c|} 
 \hline
 $k$ \textbackslash \ $n$ & 0 & 1 & 2 & 3 & 4 & 5 & 6 & 7 & 8 & 9 & 10 & 11 & 12 & 13 & 14 & 15 \\
  \hline
 0 & 1 & 1 & 1 & 1 & 1 & 1 & 1 & 1 & 1 & 1 & 1 & 1 & 1 & 1 & 1 & 1 \\ 
 \hline
 1 & &1 & 0 & 1 & 0 & 1 & 0 & 1 & 0 & 1 & 0 & 1 & 0 & 1 & 0 & 1 \\ 
 \hline
 2 & & & 1 & 1 & 2 & 2 & 3 & 3 & 4 & 4 & 5 & 5 & 6 & 6 & 7 & 7 \\ 
 \hline
 3 & & &  & 1 & 0 & 2 & 0
 & 3 & 0 & 4 & 0 & 5 & 0 & 6 & 0 & 7 \\ 
 \hline
 4 & & &  &  & 1 & 1 & 3 & 3 & 6 & 6 & 10 & 10 & 15 & 15 & 21 & 21 \\ 
 \hline
 5 & & &  &  &  & 1 & 0 & 3 & 0 & 6 & 0 & 10 & 0 & 15 & 0 & 21 \\ 
 \hline
 6 & & &  &  &  &  & 1 & 1 & 4 & 4 & 10 & 10 & 20 & 20 & 35 & 35 \\ 
 \hline
 7 & & &  &  &  &  &  & 1 & 0 & 4 & 0 & 10 & 0 & 20 & 0 & 35 \\ 
 \hline
 8 & & &  &  &  &  &  &  & 1 & 1 & 5 & 5 & 15 & 15 & 35 & 35 \\ 
 \hline
 9 & & &  &  &  &  &  &  &  & 1 & 0 & 5 & 0 & 15 & 0 & 35 \\ 
 \hline
 10 & & &  &  &  &  &  &  &  &  & 1 & 1 & 6 & 6 & 21 & 21\\  
 \hline
\end{tabular}
\bigskip
\caption{Values of $\sigma(k,n)$}
\end{table}

\begin{thm}
\label{main1}
Let $1\leq k \leq n$. Then
$$
[\Gr_k(\C^n)] \cdot [\Gr_k(\C^n)]=\sigma(k,n)
$$
while $[\Gr_k(\C^n)] \cdot \tau_*[\Gr_k(\C^n)]=0$
\end{thm}
\begin{proof}
The second statement follows immediately from the fact that 
$$ \Gr_k(\C^n) \cap \tau( \Gr_k(\C^n))=\emptyset. $$

For $1\leq k \leq n$, let 
$$s(k,n)=[\Gr_k(\C^n)] \cdot [\Gr_k(\C^n)] \in H_0(\Gr_{2k}^+(\R^{2n}))=\Z.$$
This integer can be computed by 
evaluating the Euler class of the normal bundle to 
$\Gr_k(\C^n)$ in $\Gr_{2k}^+(\R^{2n})$ 
on the fundamental class $[\Gr_k(\C^n)]$. We set $s(0,n)=1$ for all $n\geq 0$
and $s(k,n)=0$ if $k>n$ or $k<0$.

For $P \in \Gr_k(\C^n)$ the canonical isomorphism 
$$T_P \Gr_{2k}^+(\R^{2n}) = \Hom(P,P^\perp)$$
identifies the subspace $T_P\Gr_k(\C^n)$ with $\Hom_\C(P,P^\perp)$.
The fiber of the normal bundle to $\Gr_k(\C^n)$ 
at $P$ is canonically identified (as an oriented vector space
cf. \eqref{orinor}) with the space 
$$
\Hom_{\overline \C}(P,P^\perp) \cong \overline{P}^* \otimes_\C P^\perp
$$
of conjugate linear maps from $P$ to $P^\perp$. 

As the hermitian metric on $\C^n$ gives a canonical isomorphism $\overline{P}^*=P$,
letting $T$\footnote{In order to simplify the notation we omit the indices $k,n$ which are 
hopefully clear from the context.} denote the tautological bundle over $\Gr_k(\C^n)$, we 
need to compute the Euler class of the vector bundle $T \otimes T^{\perp}$. Let
\begin{equation}
\label{inclusions}
i \colon \Gr_k(\C^{n-1}) \hookrightarrow \Gr_k(\C^n)
\quad \quad
j\colon \Gr_{k-1} (\C^{n-1}) \hookrightarrow \Gr_k(\C^n) 
\end{equation}
denote the inclusions given by the expressions $i(P)=P$ and $j(P)=P\oplus \C$
(when $k=1$, $j$ includes a point as the line given by the last coordinate axis). We have 
$$
i^*(T\otimes T^\perp) = T \otimes T^\perp \oplus T,
\quad \quad
j^*(T \otimes T^\perp) = T \otimes T^\perp \oplus T^\perp.
$$

Moreover, the normal bundles to $i$ and $j$ are $T^*$ and $T^\perp$ 
respectively. 

Consider the $S^1$-action on $\Gr_k(\C^n)$ determined by the action
$$ e^{i\theta} \cdot (z_1, \ldots, z_n) = (z_1, \ldots, z_{n-1}, e^{i\theta}z_n)$$
on $\C^n$.

The fixed points of this action on $\Gr_k(\C^n)$ are precisely the 
(disjoint) images of the embeddings $i$ and $j$. Let $e^{S^1}(E)$ denote
the $S^1$-equivariant Euler class of an $S^1$-equivariant vector bundle $E\to M$ 
(meaning the Euler class of the bundle $E_{hS^1} \to M_{hS^1}$ obtained by applying
the Borel construction to the $S^1$ action on $E\to M$).

By the localization formula due to Atiyah--Bott and Berline--Vergne
\cite[Theorem 30.2]{Tu}, we have 
$$
\left \langle e^{S^1}(T \otimes T^\perp) , [\Gr_k(\C^n)] \right \rangle = 
\left \langle \frac {e^{S^1}(T \otimes T^\perp \oplus T)}{ e^{S^1}(T^*)} , [\Gr_k(\C^{n-1})] \right \rangle
+ 
\left \langle \frac{e^{S^1}(T \otimes T^\perp \oplus T^\perp)}{ e^{S^1}(T^\perp)} , 
[\Gr_{k-1}(\C^{n-1})] \right \rangle.
$$
Note that, as the degrees of the equivariant cohomology 
classes being integrated are equal to 
the dimension of the manifold on which they are being integrated, the evaluations amount
to usual evaluations of non-equivariant classes
(which are the constant terms in the equivariant classes).

Since the Euler class is multiplicative with respect to Whitney sum 
and $e^{S^1}(T^*) = (-1)^k e^{S^1}(T)$ (the Euler class is the 
top Chern class), the previous formula establishes the recursion formula
\begin{equation}
\label{recursion}
s(k,n)= (-1)^k s(k,n-1)  + s(k-1,n-1)
\end{equation}
for the intersection numbers. Note that \eqref{recursion} holds even when $k=1$ (and 
even when, in addition, $n=1$) as, in that case, the bundle $T\otimes T^\perp$ 
over the isolated fixed point $0\oplus \C$ is $0$ so that the second term in the
Atiyah-Bott-Berline-Vergne formula evaluates to $1=s(0,n-1)$. In fact,
equation \eqref{recursion} also holds for $k=0$ with our convention that $s(j,n)=0$ for $j<0$.

The formula \eqref{recursion} inductively determines the values of $s(k,n)$ from the 
values $s(0,1)=s(1,1)=1$ and one easily checks that this implies that 
$$
s(k,n)=\sigma(k,n)
$$
as required.
\end{proof}

We will now use the previous Theorem to compute the remaining homological 
intersection numbers.

\begin{prop}
\label{oppor}
Let $J' \in \JJ(\R^{2n})$ be a complex structure not compatible with the 
orientation on $\R^{2n}$ and let $1\leq k \leq n$. Then
$$
[\Gr_k(\C^n)] \cdot [\Gr^{J'}_{2k}(\R^{2n})] = \begin{cases}
\sigma(k,n) & \text{ if } n \text{ is odd and } k \text{ is even} \\
0 & \text{ if } n \text{ is odd and } k \text{ is odd} \\
\sigma(k,n-1) & \text{ if } n \text{ is even}
\end{cases} \quad \quad  = \sigma(k,n-1)
$$
and
$$
[\Gr_k(\C^n)] \cdot \tau_* [\Gr^{J'}_{2k}(\R^{2n})] = \begin{cases}
0 & \text{ if } n \text{ is odd and } k \text{ is even} \\
\sigma(k,n) & \text{ if } n \text{ is odd and } k \text{ is odd} \\
\sigma(k-1,n-1) & \text{ if } n \text{ is even}
\end{cases} \quad \quad = \sigma(k-1,n-1). 
$$
\end{prop}
\begin{proof}
Let us consider first the case when $n$ is odd. Then we can compute the 
required intersection by picking $J_0'$ in the same path component as the standard 
complex structure $J_0$ so that 
$\Gr_{2k}^{J_0'}(\R^{2n})$  and $\Gr_{2k}^{J_0}(\R^{2n})$ intersect transversally 
and choosing $J'=-J_0'$. The planes stabilised by $J'$ and $J_0'$ are the
same and $J'$ will induce the same or opposite orientation as $J_0'$ on 
a $2k$-plane $P$ according to whether $k$ is even or odd, i.e.
$$
\Gr_{2k}^{J'}(\R^{2n})=\begin{cases} 
\Gr_{2k}^{J_0'}(\R^{2n}) & \text{ if } k \text{ is even,}\\
\tau(\Gr_{2k}^{J_0'}(\R^{2n})) & \text{ if } k \text{ is odd.}
\end{cases}
$$
The orientations on the complex Grassmannians at a plane $P$ are determined by the complex
structures on $\Hom^J_\C(P,P^\perp)$, which, because $n$ is odd, is an even 
dimensional complex vector space. When $k$ is even, $J'=-J_0'$ and $J_0'$ induce the same 
orientation on $\Hom^{J_0'}_\C(P,P^\perp)=\Hom^{J'}_\C(P,P^\perp)$
hence 
$[\Gr_{2k}^{J'}(\R^{2n})]=[\Gr_{2k}^{J_0'}(\R^{2n})] \in H_*(\Gr_{2k}^+(\R^{2n})).$
When $k$ is odd, the derivative of $\tau$ at $P\in \Gr_{2k}^{J_0'}(\R^{2n})$ is given
in local coordinates by the identity map
$$
 \Hom_{\C}^{J_0'}(P,P^\perp) \xrightarrow{\id} \Hom_{\C}^{J'}(P,P^\perp)
$$
which is a conjugate linear map between these two even-dimensional complex vector spaces.
Therefore $\tau$ preserves orientation and we conclude that 
$$
[\Gr_{2k}^{J'}(\R^{2n})] =\tau_* [\Gr_{2k}^{J_0'}(\R^{2n})].
$$
The claimed result thus follows from Theorem \ref{main1}.

Now assume $n$ is even and set $J'=J_{0|\R^{2n-2}} \oplus (-J_{0|\R^2})$ to be 
the complex structure that acts as the standard complex structure
$J_0$ on the first $(n-1)$ complex coordinates and as 
$-J_0$ on the last complex coordinate.
Let $P \in \Gr_k(\C^n) \cap \Gr_{2k}^{J'}(\R^{2n})$. 
If $v \in P$ then $\frac 1 2(v- J_0 J' v) \in P$ and this expression
computes the projection of $v$ onto $\R^{2n-2} \subset \R^{2n}$.
Therefore, either $P \subset \R^{2n-2}$ or $P=P'\oplus \R^2$ for $P'$ a $(2k-2)$-plane
in $\R^{2n-2}$. Hence 
$$
\Gr_k(\C^n) \cap \Gr_{2k}^{J'}(\R^{2n}) = i(\Gr_k(\C^{n-1})) \coprod
j(\Gr_{k-1}(\C^{n-1}))$$
where $i$ and $j$ are the maps described in \eqref{inclusions}.
Moreover this decomposition corresponds to the decomposition of the intersection
into planes in which the orientations induced by $J_0$ and $J'$ agree and planes
where they do not.

By the proof of Proposition \ref{generic}, cf. also Remark \ref{genericrem}, we
can pick $g$ in $SO(2n-2)$ such that $g \Gr_{k-1}(\C^{n-1})$ intersects 
$\Gr_{k-1}(\C^{n-1})$ transversally and $g \Gr_k(\C^{n-1})$
also intersects $\Gr_{k}(\C^{n-1})$ transversally. Moreover, we can 
pick $g$ arbitrarily close to the identity and then the orientations
induced on the planes of intersection and their orthogonal complements
by the complex structures $gJ_0g^{-1}$ and $J_0$ 
will agree (as $g\Gr_l(\C^{n-1})$ will not intersect $\tau \Gr_l(\C^{n-1})$, and
sufficiently close complex structures on a plane induce the same orientation).
Let 
$$ S=\{P_1,\ldots,P_m\} =g \Gr_{k-1}(\C^{n-1}) \cap \Gr_{k-1}(\C^{n-1}), $$
$$ T=\{Q_1, \ldots, Q_r\} = g\Gr_{k}(\C^{n-1})\cap \Gr_{k}(\C^{n-1}). $$
We will show that, for $J=gJ_0g^{-1}$ (where $g$ is regarded as an element of $SO(2n)$
via the obvious inclusion $SO(2n-2)\subset SO(2n)$)
$$
 \Gr_{2k}^J(\R^{2n}) \cap \Gr_{2k}^{J'}(\R^{2n}) = i(T) \coprod j(S),
$$
that the intersections at these points are transverse and have the same
local intersection index as the corresponding points in the Grassmannians
of planes in $\R^{2n-2}$. This will complete the proof.

Any element $P \in \Gr_{2k}^{J'}(\R^{2n})$ not in the image of the inclusion 
$j \colon \Gr_{2k-2}^+(\R^{2n-2}) \to \Gr_{2k}^+(\R^{2n})$ 
projects to a $2k$-dimensional plane $P'\subset \R^{2n-2}$, 
and is therefore the graph of a unique linear
map $A \colon P' \to \R^2$. We have $P=\grph(A)\in \Gr_{2k}^J(\R^{2n}) \cap
\Gr_{2k}^{J'}(\R^{2n})$ if and only if 
the map $A$ is complex linear simultaneously as a map 
$$
(P', gJ_0g^{-1}) \to (\R^2, J_0) \quad \text{ and } \quad 
(P', J_0) \to (\R^2, -J_0).
$$
For $g$ sufficiently close to $\Id$, the subspace $\Hom_{\C}^{gJ_0g^{-1}}(P',\R^2)
\subset \Hom(P',\R^2)$
will be close to $\Hom_{\C}^{J_0}(P',\R^2)$ and will therefore intersect the plane 
$\Hom_{\C}^{(J_0,-J_0)}(P',\R^2)$ trivially. We conclude that $A$ must be the zero 
map\footnote{This also follows from the classification of pairs of orthogonal complex 
structures in the next section.}, so that $P=\grph(A)$ is in the image of 
$i\colon \Gr_{2k}^+(\R^{2n-2}) \hookrightarrow
\Gr_{2k}^+(\R^{2n})$. We conclude that $\Gr_{2k}^J(\R^{2n}) \cap
\Gr_{2k}^{J'}(\R^{2n})$ is contained in $\im(i)\coprod \im(j)$ and therefore
consists precisely of $i(T)\coprod j(S)$.

At a plane $P_i \in S$, the orthogonal complement to $j(P_i)$ in $\R^{2n}$
equals the orthogonal complement of $P_i$ in $\R^{2n-2}$, and thus
\begin{equation}
\label{decomptg}
T_{j(P_i)} \Gr_{2k}(\R^{2n}) = T_{P_i} \Gr_{2k-2}(\R^{2n-2}) \oplus 
\Hom(\R^2,P_i^\perp).
\end{equation}
By assumption, the planes $T_{P_i} \Gr_{2k-2}^J(\R^{2n-2})$ and 
$T_{P_i} \Gr_{2k-2}^{J'}(\R^{2n-2})$ add up to $T_{P_i} \Gr_{2k-2}^+(\R^{2n-2})$.
As in the previous paragraph, the intersection 
$\Hom_\C^{(J_0,gJ_0g^{-1})}(\R^2,P_i^\perp) \cap \Hom_\C^{(-J_0,J_0)}(\R^2,P_i^\perp)$
(arising from decompositions of the tangent spaces to the 
complex Grassmannians analogous to \eqref{decomptg})
is trivial so these two planes add up to $\Hom(\R^2,P_i^\perp)$ showing 
that the complex Grassmannians intersect transversely at $j(P_i)$.

Finally, let $(\phi_1, \ldots, \phi_{2(k-1)(n-k)}, \psi_1, \ldots \psi_{2(k-1)(n-k)})$ 
be a basis for 
$$T_{P_i} \Gr_{2k-2}^J(\R^{2n-2}) \oplus T_{P_i} \Gr_{2k-2}^{J'}(\R^{2n-2})  $$
compatible with the direct sum and the 
orientation determined by the complex structures on the
summands. Then a basis with the same properties for 
$$T_{j(P_i)} \Gr_{2k}^{J}(\R^{2n}) \oplus T_{j(P_i)} \Gr_{2k}^{J'}(\R^{2n})  $$
is given by 
$$(\phi_1, \ldots, \phi_{2(k-1)(n-k)},\alpha_1, \ldots, \alpha_{2n-2k},\psi_1, \ldots \psi_{2(k-1)(n-k)},\beta_1, \ldots,\beta_{2n-2k})$$
where the $\alpha_i$ and $\beta_i$ are real bases for $\Hom_\C^{(J_0,gJ_0g^{-1})}(\R^2,P_i^\perp)$ and $\Hom_\C^{(-J_0,J_0)}(\R^2,P_i^\perp)$ respectively. 
Since $g$ is close to the identity, the argument near \eqref{orinor} implies that 
the basis $(\alpha_1, \ldots, \alpha_{2n-2k},\beta_1,\ldots,\beta_{2n-2k})$ has the 
same orientation as the standard basis for $\Hom(\R^2, P^\perp)$. This implies
that the local intersection indices at $P_i$ and $j(P_i)$ agree.

The argument for the points in $i(T)$ is entirely analogous and is therefore omitted.
\end{proof}

\begin{remark}
Theorem \ref{main1} and Proposition \ref{oppor} can alternatively be proved by  
picking explicit pairs of orthogonal complex structures whose corresponding complex 
Grassmannians intersect transversally (any generic pair as in the statements of 
Corollaries \ref{genorthopos} and \ref{genorthoneg} will do)
and then determining the local intersection numbers at the points of intersection
of the Grassmannians.
\end{remark}

From Proposition \ref{generic}, Theorem \ref{main1}, and Proposition \ref{oppor},
we obtain the following result.
\begin{cor}
\label{estints}
Let $V$ be a real vector space of dimension $2n$. There exists an open dense set of pairs 
$(J_0,J_1)$ in $\mathcal J(V)^2$ for which there are at least 
$\sigma(k,n)$ $2k$-planes in $V$ which are stabilised by both $J_0$ and $J_1$ when
these induce the same orientation on $V$ or $n$ is odd, and at least 
$\sigma(k,n-1)+\sigma(k-1,n-1)$ stabilised $2k$-planes when $J_0$ and $J_1$ induce
opposite orientations and $n$ is even.
\end{cor}
Note that when $n$ and $k$ are both even $\sigma(k,n-1)+\sigma(k-1,n-1)=\sigma(k,n)$.

\section{The classification of pairs of orthogonal complex structures}

Two pairs of complex structures $(J_0,J_1)\in \JJ(V), (J_0',J_1') \in \JJ(V')$ 
are isomorphic if 
there exists a linear isomorphism $\phi \colon V \to V'$ such that 
$J_i'=\phi J_i \phi^{-1}$ for $i=0,1$. If the vector spaces are Euclidean
and the complex structures orthogonal we further require that $\phi$ be an orthogonal map.

In the next section we will completely describe the possible intersections
of pairs of orthogonal complex structures. This will rely on the classification
of pairs of orthogonal complex structures on a Euclidean 
even dimensional real vector space $V$ up to isomorphism, which is our goal 
in this section. We will write $J(V)$ for 
the space of orthogonal complex structures on $V$.

\subsection{Dimensions two and four}
On $\R^2$ there are exactly two orthogonal complex structures, corresponding to 
the two possible orientations. Identifying $\R^2$ with $\C$, these are left 
multiplication by $i$ and by $-i$. It follows that, up to isomorphism, 
there are exactly two pairs of 
complex structures, namely $(i,i)$ and $(i,-i)$ on the vector space $\C$ 
(note that complex conjugation gives an 
isomorphism between $i$ and $-i$). We will abbreviate these two
isomorphism classes of pairs by $\C$ and $\overline \C$ respectively.

On $\R^4$, the orthogonal complex structures can be identified with left and right 
multiplication by unit imaginary quaternions. To be specific, we identify $\R^4$ with 
$\C^2$ in the standard way and the latter with $\HH$ via
$$
(z_1, z_2) \mapsto z_1 + z_2 j
$$
Left multiplication by a unit imaginary quaternion is compatible with the standard 
orientation on $\C^2$, while right multiplication is not. 

\begin{defn}
Let $0<\theta<\pi$. Then
$$\HH_\theta = (\HH, i,i e^{j\theta})$$ 
denotes the Euclidean vector space $\HH$ together with the pair of orthogonal 
complex structures given by left multiplication by $i$ and left multiplication 
by $i e^{j\theta}= \cos(\theta) i + \sin(\theta) k$.
\end{defn}

Below, we will abbreviate $(\HH, i, i e^{j\theta})$ by 
$(i,i e^{j\theta})$ leaving the vector space $\HH$ implicit.

The above definition makes sense for any real $\theta$ (which we may assume to lie in 
$[0,2\pi]$). However $\HH_0=\C \oplus \C$, $\HH_\pi= \overline \C \oplus \overline\C$ 
and we will see shortly that $\HH_\theta\cong \HH_{2\pi-\theta}$.

Note that $\HH_\theta \not \cong \HH_{\theta'}$ for $0< \theta \neq \theta'< \pi$
as the eigenvalues of the orthogonal endomorphism $-J_0J_1 = e^{j\theta}$, represented
in the standard basis by the matrix
$$
\begin{bmatrix}
\cos \theta & 0 & -\sin \theta & 0 \\
0 & \cos \theta & 0 & \sin \theta \\
\sin \theta & 0 & \cos \theta & 0 \\
0 & -\sin \theta & 0 & \cos \theta
\end{bmatrix}
$$
are $\pm e^{i\theta}$ (both with multiplicity $2$).

\begin{prop}
\label{dim4}
The set
$$ \{ \C^2, \overline \C^2, \C \oplus \overline \C\} \cup \{ \HH_\theta \colon 
0<\theta<\pi \} $$
contains exactly one representative of each isomorphism class of pairs of
orthogonal complex structures on Euclidean vector spaces of dimension $4$.
\end{prop}

\begin{remark}
In particular, there is a unique isomorphism class of pairs of orthogonal 
complex structures on $\R^4$ inducing opposite orientations, 
namely $\C \oplus \overline \C$.
\end{remark}

\begin{proof}[Proof of Proposition \ref{dim4}]
The eigenvalues of $-J_0J_1$ distinguish all pairs of complex structures
listed. In order to check that the isomorphism classes are exhausted 
by the ones in the statement we may assume that the underlying vector space is $\HH$. 

Consider first the case where $J_0$ and $J_1$ induce the same orientation.
We may assume that the induced orientation is the standard one and then 
an isomorphism between two such pairs will necessarily be an element
$g\in SO(4)$, which can be written in the form 
$$ x \mapsto q_1 x \overline{q_2} $$
with $q_1,q_2$ unit quaternions. We may further
assume that $J_0$ is left multiplication by $i$ and $J_1$ is left multiplication
by the purely imaginary unit quaternion $\alpha$. Two such pairs $(i,\alpha)$ and
$(i,\alpha')$ are isomorphic if and only if there exist unit quaternions $q_1$ 
and $q_2$ such that 
$$
q_1 i x \overline {q_2} = iq_1 x \overline{q_2} \quad \text{ and }\quad 
q_1 \alpha x \overline {q_2} =  \alpha' q_1 x \overline{q_2}
$$
for all $x \in \HH$. The first condition amounts to  
$q_1= \cos \theta + i \sin \theta$ for some $\theta$ while the second 
states that $\alpha' =(\cos \theta + i \sin \theta)\alpha(\cos \theta -i \sin \theta)$.
Writing $\alpha= i \cos \phi + z_2 j$ with $0\leq \phi \leq \pi$ and $|z_2|=|\sin \phi|$
the second condition becomes
$$
\alpha'= i \cos \phi + e^{i2\theta} z_2 j
$$
It follows that two pairs $(i,\alpha)$ and $(i,\alpha')$ are isomorphic if and 
only if the components of $\alpha$ and $\alpha'$ along $i$ are equal. This completes the 
proof in the equal orientation case. 

If $J_0$ and $J_1$ have opposite orientations then we can again assume that $J_0$
is left multiplication by $i$. $J_1$ must then be given by right multiplication by 
a unit imaginary quaternion $\alpha$. Writing $\alpha= i \cos \beta + w j$, the equation
$$
-i(z_1+z_2j) = (z_1+z_2j)(i \cos \beta + w j)
$$
translates into the homogeneous system for $(z_1,z_2)$ with matrix 
$$
\begin{bmatrix}
-i(1+\cos \beta) & \overline w \\
-w & -i(1-\cos \beta)
\end{bmatrix}
$$
which has zero determinant. Hence there is a complex line in $\HH$ where right
multiplication by $\alpha$ agrees\footnote{This is also a special case of Lemma 
\ref{line}.} with left multiplication by $(-i)$. On the orthogonal line, right
multiplication by $\alpha$ must equal left multiplication by $i$ or $-i$
and the incompatibility with the orientation means it must be $i$.
\end{proof}

\subsection{Classification}

The following point of view seems useful in the study of pairs of orthogonal almost 
complex structures: consider the amalgam
$$G=\Z/4 \ast_{\Z/2} \Z/4.$$
A pair of orthogonal complex structures $(J_0,J_1)$ on $V$ can be 
identified with an orthogonal (real) representation of $G$ on $V$ 
where the generator of the canonical subgroup $\Z/2 \subset G$ 
acts as $-\Id$; namely, a pair $(J_0,J_1)$ is identified with the representation 
which assigns $J_0$ to the 
first canonical generator $1\in \Z/4$ and $J_1$ to the second. Conversely,
such a representation of $G$ is determined by the isomorphisms which are assigned 
to the canonical generators of $\Z/4$ and these must be orthogonal complex structures.

\begin{defn}
An orthogonal representation of $G=\Z/4 \ast_{\Z/2} \Z/4$ is said to be 
\emph{admissible} if it restricts along the canonical map $\Z/2 \to G$ to a representation where the generator acts by $-\Id$.
\end{defn}
The admissible representations $V$ of $G$ break up into two 
classes, according to whether the corresponding complex structures 
$J_0$ and $J_1$ induce the same, or opposite orientations of $V$.
We will refer to the former as \emph{positive} admissible representations 
and to the latter as \emph{negative}. Note that the direct sum of two positive
or negative admissible representations is positive, while the direct sum of a positive
with a negative representation is negative.

We remark that the admissible representations $\C,\overline \C$ and $\HH_\theta$ (with
$0<\theta<\pi$) defined in the previous section are all irreducible. With the exception
of $\overline \C$ they are positive. 
The next proposition shows that all admissible representations are direct sums of these.

\begin{lemma}
\label{line}
Let $J_0,J_1$ be orthogonal complex structures on the Euclidean vector space $V$
inducing opposite orientations. Then there exists a plane $L \subset V$
stabilised by both $J_0$ and $J_1$ on which $J_1=-J_0$.
\end{lemma}
\begin{proof}
By Proposition \ref{oppor}, $[\Gr_1^{J_1}(V)]\cdot \tau_*[\Gr_1^{J_0}(V)]=1$ for any
$V$. By Proposition \ref{generic} this implies that for an open dense set of 
$J$'s in the component of $J(V)$ not containing $J_0$, the endomorphism $J+J_0$ vanishes
on at least one complex line. Since the condition $\dim \ker(J+J_0) \geq 1$
is closed, the desired statement follows.
\end{proof}

\begin{remark}
Given $J_0 \in J(V)$, the sets $S_k=\{J \in J(V) \colon \dim_\C \ker(J+J_0)= k\}$ form a 
stratification of $J(V)$. One can show that $S_1$ is the open dense stratum in the 
component of $J(V)$ which does not contain $J_0$ and this gives an alternative proof
of Lemma \ref{line}.
\end{remark}

\begin{prop}
\label{class}
Every admissible representation of $G$ is isomorphic to a direct sum of the admissible, 
irreducible representations of dimensions $2$ and $4$ (which are $\C,\overline \C$
and $\HH_\theta$ for $0<\theta<\pi$).
\end{prop}
\begin{proof}
We consider the case of positive admissible representations first. 

The space of pairs of orthogonal complex structures  
inducing a given orientation on a vector space of dimension $2n$ can be identified with
the homogeneous space $(SO(2n)/U(n))^2$. 
Two pairs are isomorphic if and only if they are in the same diagonal $SO(2n)$ orbit.
Thus the space of isomorphism classes of pairs is
$$
SO(2n)\backslash((SO(2n)/U(n))^2) = U(n)\backslash SO(2n)/U(n)
$$
Letting $J_0$ denote the standard complex structure on $\R^{2n}$,
the orbit of $g$ in $SO(2n)$ in the biquotient corresponds to the pair of complex 
structures $(J_0, gJ_0g^{-1})$.

Let $\lp$ denote the orthogonal complement to $\uu(n) \subset \so(2n)$ (i.e. 
the space of skew-symmetric, complex anti-linear matrices) and 
$\la\subset \lp$ be a maximal abelian Lie algebra (amongst subalgebras of $\so(2n)$ 
contained in $\lp$).
By \cite[Theorem VII.8.6]{He} (see also \cite[Section 9]{E}), the map
$$\la \hookrightarrow \so(2n)/\uu(n) \xrightarrow{\exp} SO(2n)/U(n) \to 
U(n)\backslash SO(2n)/U(n)$$
is surjective. 

An explicit choice for $\la$ is the set of matrices 
$$
\begin{bsmallmatrix}
0&0  & & & & & & & a_1 & 0 \\
0&0  & & & & & & & 0 & -a_1\\
& & \ddots & & &  & & \iddots & & \\
& &  & 0&0 & a_{\frac n 2} & 0 & &   &  \\
& &  & 0&0 & 0 & -a_{\frac n 2} & & & \\
& & &  -a_{\frac n 2} & 0 & 0 & 0& &   & \\
& & &  0& a_{\frac n 2} & 0& 0& & &   \\
& & \iddots & & & & &\ddots  & &\\
-a_1 & 0 & & & & & & & 0& 0\\
0 & a_1 & & & & & & &0 & 0
\end{bsmallmatrix}
$$
when $n$ is even and 
$$
\begin{bsmallmatrix}
0&0 & & & & & & & a_1 & 0 \\
0&0 & & & & & & & 0 & -a_1\\
& & \ddots& & & &  & \iddots & & \\
& & & & & & a_{\frac{n-1} 2} & 0 &  &  \\
& & & & &  & 0 & -a_{\frac {n-1} 2} & & \\
& & & & 0& 0& 0 &  &  &  \\
& & & & 0 & 0 &  &  & & \\
& & -a_{\frac {n-1} 2} & 0 &  & &  & & & \\
& &  0& a_{\frac{n-1} 2} & & & &  & & \\
& &  \iddots &  & & & &  \ddots & & \\
-a_1 & 0 & & & & & & & 0& 0\\
0 & a_1 & & & & & & &0 & 0
\end{bsmallmatrix}
$$
when $n$ is odd (with $a_i \in \R$). The subalgebra $\la$ 
has dimension $\lfloor \frac n 2 \rfloor$ - the rank of the symmetric space 
$SO(2n)/U(n)$.

Thus $\la$ consists of complex anti-linear, anti-diagonal, off-diagonal matrices 
which (when the $2\times 2$ blocks are regarded as complex numbers in the natural way)
have real entries. If we write $\C^{n}\cong \HH^{\frac n 2}$ for $n$ even,
respectively $\C^n \cong \HH^{\frac{n-1}{2}} \oplus \C$ for $n$ odd, 
by identifying
$$
(z_1, \ldots,z_n) \leftrightarrow (z_1+z_n j, z_2+z_{n-1}j, \ldots, z_{\frac n 2} + z_{\frac n 2 + 1} j  )
$$
when $n$ is even and 
$$
(z_1, \ldots,z_n) \leftrightarrow (z_1+z_n j, z_2+z_{n-1}j, \ldots, z_{\frac{n-1} 2} + z_{\frac{n+3} 2 }j, z_{\frac{n+1}{2}} )
$$
when $n$ is odd, the above matrices represent left multiplication by 
the diagonal matrices in $\End(\HH^{\frac n 2})$, respectively 
$\End(\HH^{\frac{n-1} 2} \oplus \C)$, 
$$
\begin{bmatrix}
-a_1 j & & \\
& \ddots & \\
& & -a_{\frac n 2} j
\end{bmatrix}
\quad 
\text{( $n$ even)} \quad \quad \quad 
\begin{bmatrix}
-a_1 j & & & \\
& \ddots & & \\
& & -a_{\frac{n-1}{2}} j & \\
& & & 0
\end{bmatrix}
\quad \text{ ($n$ odd).}
$$
The exponentials of these matrices are
$$
\begin{bmatrix}
e^{-a_1 j} & & \\
& \ddots & \\
& & e^{-a_{\frac n 2} j}
\end{bmatrix}
\quad 
\text{( $n$ even)} \quad \quad \quad 
\begin{bmatrix}
e^{-a_1 j} & & &\\
& \ddots & & \\
& & e^{-a_{\frac{n-1} 2} j}& \\
& & & 1
\end{bmatrix}
\quad \text{ ($n$ odd).}
$$
The space $J(\R^{2n})$ is identified with $O(2n)/U(n)$ by acting with an element of 
$O(2n)$ on the 
standard complex structure $J_0$ via conjugation. Therefore, when regarded as
complex structures on $\R^{2n}=\HH^{\lfloor \frac n 2\rfloor} (\oplus \C)$, 
the images of the previous matrices in $SO(2n)/U(n)$ are given
by coordinatewise left multiplication by the imaginary unit quaternions
$$
e^{-a_l j} i e^{a_l j} = i e^{2 a_l j} = \cos(2 a_l) i + \sin(2 a_l) k.
$$
We conclude that the image of $\la$ in $U(n)\backslash SO(2n)/U(n)$ 
corresponds to the set of pairs of complex structures on 
$\HH^{\lfloor \frac n 2 \rfloor}(\oplus \C)$ given by the expression
\begin{equation}
\label{expJ}
(i, ie^{2a_1 j}) \oplus \cdots \oplus (i,ie^{2a_{\lfloor \frac{n}{2}\rfloor} j}) 
\left( \oplus (i,i) \right) )
\end{equation}

Note (cf. proof of Proposition \ref{dim4}) that for $\pi < \theta < 2\pi$ we have 
$(i,ie^{\theta j}) \cong (i, ie^{(2\pi-\theta)j}) \cong \HH_{2\pi - \theta}$ 
as $ie^{\theta j}$ and $ie^{(2\pi-\theta)j}$ have the same component along $i$. 
For $\theta=0$ or $\pi$ we
have $(i,ie^{\theta j})=(i,i)\cong \C^2$ or $(i,-i)\cong \overline \C^2$ respectively. 
We conclude that all positive admissible representations
may be expressed as sums of copies of $\C$, $\overline\C^2$ or $\HH_\theta$ with 
$0<\theta<\pi$ as required. 

Suppose now that $V$ is a negative admissible representation corresponding to a
pair $(J_0,J_1)$ of orthogonal complex structures. As $J_1$ is not compatible with the orientation induced by $J_0$, by Lemma \ref{line}, 
$J_1=-J_0$ on some line. Therefore we may decompose $V$ as an orthogonal direct sum 
$\overline \C \oplus W$ with $W$ positive admissible. This completes the proof.
\end{proof}

\begin{remark}
In light of Proposition \ref{class} we see that the isomorphism type of a pair $(J_0,J_1)$
of orthogonal complex structures is completely determined by the eigenvalues of the
orthogonal map $-J_0J_1$ and their multiplicities (cf. beginning of the proof of Proposition \ref{dim4}).
\end{remark}

\section{Intersections of orthogonal complex structures}
Armed with the classification of the previous section we will now describe in 
detail the intersection of two complex Grassmannians $\Gr_{2k}^{J_i}(\R^{2n})$
for a pair of orthogonal complex structures $(J_0,J_1)$.

As in the previous section, $G$ denotes the amalgam $\Z/4 \ast_{\Z/2} \Z/4$.
\begin{lemma}
\label{endos}
The algebra of endomorphisms of the admissible irreducible representations of 
$G$ are 
$$
\End(\C) = \C, \quad \quad \End(\overline \C)=\C, \quad \quad \End(\HH_\theta) = \HH
\quad (0<\theta<\pi).
$$
\end{lemma}
\begin{proof}
This is clear for the representations of dimension $2$. For the $4$-dimensional
representations this is a consequence of the fact that $i$ and $i \cos \theta  + 
k\sin \theta$ generate $\HH$ as an algebra over $\R$ when $0<\theta<\pi$.
\end{proof}

By Schur's lemma, the only map between non-isomorphic irreducible 
admissible representations is the zero map (note that a sub-representation or a 
quotient of admissible representations is still admissible). 

\begin{thm}
\label{descint}
Let $V$ be a real euclidean oriented vector space and $(J_0,J_1)$ be a pair
of orthogonal complex structures with isomorphism class 
$$\HH_{\theta_1}^{r_1} \oplus \cdots \oplus \HH_{\theta_m}^{r_m} \oplus \C^\ell 
\oplus \overline{\C}^s$$ 
where $0<\theta_1< \ldots <\theta_m < \pi$, $r_i>0$ and $m,\ell, s\geq 0$.

Then, for $k<n$ the space of oriented $2k$-planes $P \subset V$ which 
are stabilised by both $J_0$ and $J_1$ is
$$
\coprod_{\stackrel{2t_1+\ldots +2t_m + \ell' + s'=k}
{0\leq t_i \leq r_1, 0 \leq \ell' \leq \ell, 0 \leq s' \leq s}} \
\Gr^{\HH}_{t_1}(\HH^{r_1}) \times \cdots \times \Gr^{\HH}_{t_m}(\HH^{r_m}) \times \Gr^{\C}_{\ell'}(\C^\ell) \times \Gr^\C_{s'}(\C^s) \quad \subset
\quad 
\Gr_{2k}^+(V)
$$
The components consisting of planes where $J_0$ and $J_1$ induce the same orientation
are those for which the index $s'$ is even.
\end{thm}
\begin{proof}
A $2k$-plane $P\subset V$ which is invariant under $J_0$ and $J_1$ is isomorphic as 
an admissible $G$-representation to 
$$ \HH_{\alpha_1}^{t_1} \oplus \cdots \oplus \HH_{\alpha_k}^{t_k} \oplus
\C^{\ell'} \oplus \overline\C^{s'}$$
with the $0<\alpha_i <\pi$ distinct, $k\geq 0$, $t_i\geq 1$ and $\ell',s'\geq 0$.
By Schur's lemma the set $\{\alpha_1,\ldots,\alpha_k\}$ must be contained
in $\{\theta_1, \ldots, \theta_m\}$ and each summand of $P$ must map non-trivially
only to a summand of $V$ of the same isomorphism type. As $P\subset V$ we must 
further have for $\alpha_i=\theta_j$ that $t_i \leq r_j$.

By \Cref{endos} the maps between direct sums of copies 
of a given irreducible representation are matrix algebras over $\C$ or $\HH$. 
The subspaces of a $p$-fold direct sum of an irreducible 
which are isomorphic to a $q$-fold direct 
sum ($q\leq p$) can therefore be identified with the quotient of full 
rank $p\times q$ matrices by the action of invertible $q\times q$ matrices, i.e. with 
the Grassmannian of $q$-planes in $\C^p$ or $\HH^p$. 

It follows that a choice of 
$P\subset V$ of a given isomorphism type 
corresponds to an independent 
choice of $P_i \in \Gr_{t_i}(\HH^{r_j}_{\theta_j})$ for each $i=1,\ldots,m$,
$P' \in \Gr^{\C}_{\ell'}(\C^\ell)$, and $P'' \in \Gr^{\C}_{s'}(\C^s)$
such that $(\sum_{i=1}^m 4t_i)+2\ell' + 2s'= \dim_\R P$. 
This is the content of the description
of the space of oriented $2k$-planes given in the statement.

Finally, a representation  $\HH_{\alpha_1}^{t_1} \oplus \cdots 
\oplus \HH_{\alpha_k}^{t_k} \oplus \C^{\ell'} \oplus \overline\C^{s'}$ is positive
if and only if $s'$ is even.
\end{proof}

It is worthwhile highlighting the generic cases in the previous Theorem, which apply to
an open dense set in the space of pairs of orthogonal complex structures. These correspond
to the principal orbit types for the action of $U(n)$ on $O(2n)/U(n)$ which have
isotropy type $\Sp(1)^{\lfloor \frac n 2 \rfloor} (\times U(1)) \subset U(n)$ 
in the $SO(2n)/U(n)$ component and isotropy type $\Sp(1)^{\lfloor \frac n 2 \rfloor-1} 
\times U(1)\times U(1)$ or $\Sp(1)^{\lfloor \frac n 2 \rfloor} \times U(1)$
in the remaining component, according to whether $n$ is even or odd.
These are the principal orbit types because the corresponding isotropy groups
(the automorphism groups of the pairs of complex structures)
have the lowest possible dimension. 

Alternatively, at least in the case when the complex structures induce 
the same orientation, these generic pairs of 
complex structures are the $SO(2n)$-orbits of pairs $(i,J)$ with $J$ 
in the image under the exponential map of a component $Q_0$ of the 
complement $\la \setminus D(SO(2n),U(n))\subset \la$ of the diagram of 
$(SO(2n),U(n))$ satisfying $0 \in \overline{Q_0}$ (cf. \cite[Theorem VII.8.6]{He}). The 
remaining pairs of complex structures are orbits where $J$ is in
the image of $\partial Q_0$; these have larger dimensional isotropy and form a 
positive codimensional subspace of the space of pairs.

\begin{cor}
\label{genorthopos}
Let $V$ be a real euclidean oriented vector space of dimension $2n$.
Let $0<\theta_1<\ldots <\theta_{\lfloor \frac n 2\rfloor} < \pi$
and $(J_0,J_1)$ be a pair of orthogonal complex structures on $V$ with isomorphism type 
$$\HH_{\theta_1} \oplus \cdots \oplus \HH_{\theta_{\lfloor \frac n 2\rfloor}} \quad 
\text{ or } \quad 
\HH_{\theta_1} \oplus \cdots \oplus \HH_{\theta_{\lfloor \frac n 2\rfloor}} 
\oplus \C$$ 
according to whether $n$ is even or odd . Then $J_0$ and $J_1$ intersect on exactly $\sigma(k,n)$ $2k$-planes and induce the same orientation on these planes of intersection.
\end{cor}
\begin{proof}
An admissible sub-representation $P\subset V$ is necessarily positive and consists
of a direct sum of a subset of the summands in the decomposition of $V$ into 
irreducibles. If $\dim_\C P=k$ there are exactly $\sigma(k,n)$ possible choices.
\end{proof}

\begin{cor}
\label{genorthoneg}
Let $V$ be a real euclidean oriented vector space of dimension $2n$.
Let $0<\theta_1<\ldots <\theta_{\lfloor \frac n 2\rfloor} < \pi$
and $(J_0,J_1)$ be a pair of orthogonal complex structures with isomorphism type 
$$\HH_{\theta_1} \oplus \cdots \oplus \HH_{\theta_{ \frac n 2-1}}
\oplus \C \oplus \overline \C \quad 
\text{ or } \quad 
\HH_{\theta_1} \oplus \cdots \oplus \HH_{\theta_{ \lfloor \frac n 2\rfloor}} 
\oplus \overline \C$$ 
according to whether $n$ is even or odd . Then, we have the following:
\begin{itemize}
\item For $n$ even, $J_0$ and $J_1$ intersect on exactly $\sigma(k,n-1)$ 
$2k$-planes inducing the same orientation 
and on $\sigma(k-1,n-1)$ $2k$-planes inducing the opposite orientation,
\item For $n$ odd,  $J_0$ and $J_1$ intersect on exactly $\sigma(k,n)$ $2k$-planes, 
inducing the same orientation when $k$ is even and inducing the opposite orientation when 
$k$ is odd.
\end{itemize}
\end{cor}
\begin{proof}
The proof is entirely analogous to that of Corollary \ref{genorthopos}.
\end{proof}

Thus, for generic pairs of orthogonal complex structures, the estimates of 
Theorem \ref{main1} and Proposition \ref{oppor} give the exact number of 
intersections. This has the following consequence.
\begin{cor}
\label{positive}
Let $J_0,J_1$ be orthogonal complex structures on $\R^{2n}$. 
If $\Gr_{2k}^{J_0}(\R^{2n})$ and $\Gr_{2k}^{J_1}(\R^{2n})$ (respectively 
$\tau(\Gr_{2k}^{J_1}(\R^{2n})$) intersect transversally 
in $\Gr_{2k}^+(\R^{2n})$ then the local intersection numbers at each point 
of intersection are $+1$.
\end{cor}

By Theorem \ref{descint} the intersection is transverse precisely in the generic
cases described by Corollaries \ref{genorthopos} and \ref{genorthoneg} (otherwise 
some component of the intersection will contain a nontrivial complex or quaternionic 
Grassmannian). The previous Corollary can also be checked explicitly in those examples,
and this gives an alternative route to Theorem \ref{main1} and Proposition 
\ref{oppor}.

\section{Remarks on intersections of non-orthogonal complex structures}

In this section we discuss some differences between the intersection properties 
of orthogonal complex structures and general complex structures. 

\begin{lemma}
Let $(J_0,J_1)$ be orthogonal complex structures on $\R^{2n}$. 
If a point $P \in \Gr_{2k}^{J_0}(\R^{2n}) \cap \Gr_{2k}^{J_1}(\R^{2n})$ is an isolated
point of intersection in $\Gr_{2k}^+(\R^{2n})$ 
then the two complex Grassmannians intersect transversally at $P$.
\end{lemma}
\begin{proof}
As the $J_i$ are orthogonal, $P^\perp$ is stabilised by both $J_i$. In the local 
coordinates $\Hom(P,P^\perp)$ at $P$, the complex Grassmannians are expressed as
the planes $\Hom_\C^{J_0}(P,P^\perp)$ and $\Hom_\C^{J_1}(P,P^\perp)$ which 
are transverse if and only if their intersection  is $\{0\}$.
\end{proof}

Example \ref{exnontrans} below shows the previous statement is not true
for general pairs of 
complex structures. Note that if $J_0,J_1 \in \JJ(V)$ intersect non-transversally 
at an isolated point $P$, it is not possible to find a complement $Q$ for $P$ in $V$
which is stabilised by both $J_0$ and $J_1$ (otherwise 
the argument in the proof of the previous 
lemma would prove the intersection to be transverse).

We have seen in Proposition \ref{generic} that the Grassmannians associated to complex 
structures generically intersect transversally at a finite set of planes.
In the orthogonal case, the local intersection numbers at each of the finite 
set of intersection point are always $+1$ (Corollary \ref{positive}). Again, we will see 
in Example \ref{exnontrans} below that this need not be true for general pairs of complex
structures: $\Gr_{2k}^{J_0}(\R^{2n})$ and $\Gr_{2k}^{J_1}(\R^{2n})$ may intersect
transversally at a finite set of points having local intersection numbers of opposite 
signs. In these cases the lower bounds obtained in section 2 are strictly less than
the number of actual intersection points.

\begin{example}
\label{exnontrans}
Recall that orthogonal complex structures on $\R^4\cong \HH$ 
compatible with the orientation can be identified with the unit sphere $S^2$
inside the imaginary unit quaternions.

For each $P \in \Gr_2^+(\R^4)$ there is a unique 
$\alpha \in S^2 \subset \im(\HH)$ such that left multiplication by 
$\alpha$ is a complex structure on $P$ compatible with the orientation\footnote{ Following the answer by user Thomas to \href{https://math.stackexchange.com/questions/1828484/grassmanian-2-4-homeomorphic-to-s2-times-s2}{https://math.stackexchange.com/questions/1828484/grassmanian-2-4-homeomorphic-to-s2-times-s2}.} (if $( v_1,v_2)$ is an oriented orthonormal basis for $P$, then 
$\alpha = v_2 \overline v_1 \in \HH$). This assignment yields a fiber bundle 
$$ \Gr_2^+(\R^4) \xrightarrow{\pi} S^2 $$
with fiber $\CP^1_\alpha=\Gr_2^\alpha(\R^4)$
over $\alpha \in S^2$ (this bundle can be identified with
the holomorphic bundle $\mathbb{P}(\mathcal O(2) \oplus \mathcal O) \to \CP^1$).

An arbitrary complex structure $J$ on $\R^4$ compatible with the orientation, 
determines a subspace
$$
\CP^1_J =\Gr_2^J(\R^4) \subset \Gr_2^+(\R^4)
$$
which is homologous to a fiber and hence has homological intersection $0$ with any fiber.
It follows that $\CP^1_J$ can not possibly intersect a fiber transversely 
and positively. Let us consider the specific example 
$$
J = \begin{bmatrix} 
0 & -\frac 1 {a^2} & & \\
a^2 & 0 & & \\
& & 0 & -\frac 1 {b^2} \\
& & b^2 & 0
\end{bmatrix}
$$
where $a,b>0$. This complex structure is obtained from the standard $J_0$ via 
conjugation by the diagonal matrix 
$$
g=\begin{bmatrix}
\frac 1 a & & & \\
& a & & \\
& &\frac 1 b & \\
& & & b
\end{bmatrix}
$$
hence $\CP^1_J = g\CP^1$. Writing $\CP^1 = \{ [x+iy:1] \} \cup \{[1:0]\}$ we have 
$$
\CP^1_J = \{ \Span \{ \left( \tfrac x a, a y, \tfrac 1 
b,0\right), \left( -\tfrac y a, a x, 0,b \right) \} \colon x,y \in \R\} \cup 
\{\R^2\oplus 0\}
$$
A point in $\CP^1_J$ parametrised by $(x,y)$ lies in $\CP^1$, i.e. projects to $i$ under 
$\pi$,  if and only if we can solve the following equation for $t\in \R$ and $s>0$ (the
latter condition amounting to the orientation induced by $i$ being the correct one):
$$
i \left( \tfrac x a, a y, \tfrac 1 b,0 \right)=\left( -a y , \tfrac{x}a, 0,\tfrac{1}b 
\right) = t \left( \tfrac x a, a y, \tfrac 1 b,0 \right) + s 
\left( -\tfrac y a, a x, 0,b \right)
$$
For $a \neq b, \frac 1 b$, this equation can only be solved when $x=y=0$. Thus $\CP^1_J$ 
and $\CP^1$ intersect exactly at $[1:0]$ 
and $[0:1]$ and it is not hard to check that the intersection is transverse 
at both points with opposite signs, corresponding to the fact that the homological 
intersection is $0$.

In general, in order to find the image $\alpha = \pm\sqrt{1-u^2-v^2}i + uj + vk \in S^2$
of a point in $g\CP^1$ under the projection $\pi$, 
we can find the unique values of $u,v$ for which the following equation for $t \in \R$
and $s>0$ can be solved
$$
\left(\pm \sqrt{1-u^2-v^2}i + uj + vk\right)\left(\frac x a + ayi + \frac j b\right) = 
t\left(\frac x a + ayi +\frac j b\right) + s\left( -\frac y a + axi + bk\right).
$$
At least for values of $a,b$ close to $1$, the points in $\CP^1_J$ will map under $\pi$
to points close to $i$ so the sign in the above equation will be positive. Then
the components along $j$ and $k$ 
determine the values of $t,s$ uniquely and, plugging them in, 
the above equation becomes
$$
\begin{cases}
\left( \frac{bx^2}{a^2} + \frac{y^2}{b} +\frac 1 b\right) u + 
xy(b-\tfrac{1}{a^2b} ) v = (\frac 1 {ab^2}-a)y \sqrt{1-u^2-v^2}, \\
xy(b- \tfrac{a^2}{b})u + (\tfrac{x^2}b + a^2by^2 + \tfrac 1 b)v=
 (\frac 1 a - \frac{a}{b^2}) x \sqrt{1-u^2-v^2},
\end{cases}
$$
which can be solved using polar coordinates for $u,v$.

At least for $a,b$ close to $1$, the image of $\CP^1_J$ under the projection $\pi$ 
is a proper compact subset $K \subset S^2$
containing $i$, and if $\alpha=i\sqrt{1-u^2-v^2}+uj+vk$ is a point in the boundary of $K$,
the intersection of the fiber $\CP^1_{\alpha}$ with $\CP^1_J$ will not be transverse.

For instance, assuming $a=b$, the second equation implies $v=0$ and the system
is equivalent to 
$$
\frac{u}{\sqrt{1-u^2}}= \left( \frac 1 b^2 - b^2 \right) \frac{y}{x^2+ y^2 +1}.
$$
The maximum of the function of $(x,y)$ on the right hand side is assumed 
at a single point. Letting $u_{max}$ denote the corresponding value of $u$,
the fiber over $i\sqrt{1-u_{max}^2} + u_{max} j$ intersects $\CP^{1}_J$ at a 
single point. As the homological intersection is zero, the intersection can not
be transverse.
\end{example}

We have seen that for pairs of orthogonal complex structures (inducing the same 
orientation say) the number of planes of intersection is either equal to 
the lower bound implied by Theorem \ref{mainno} or infinite and we have shown by example
this is no longer true for general pairs of complex structures. It seems likely
that the lower bound of Theorem \ref{mainno} holds in complete generality for 
arbitrary pairs of complex structures, not just generically, but we have not been 
able to prove this.


\begin{thebibliography}{FGM}
\bibitem[Ba]{Ba} F. Battaglia, \emph{Linear almost complex structures on $S^6$.} Ann. Mat. Pura Appl. (4) 156, 181-193, 1990.
\bibitem[E]{E} J. H Eschenburg, \emph{Lecture notes on symmetric spaces}, available at
\url{https://myweb.rz.uni-augsburg.de/~eschenbu/symspace.pdf}
\bibitem[FGM]{FGM} B. Ferlengez, G. Granja and A. Milivojevic, On the topology of the space of almost complex structures on the six sphere, New York J. of Math. 27 (2021) 1258-1273.
\bibitem[He]{He} S. Helgason, \emph{Differential Geometry, Lie groups and Symmetric spaces}, Grad. Stud. Math., 34, American Mathematical Society, Providence, RI, 2001, xxvi+641 pp.
\bibitem[Tu]{Tu} L. Tu, \emph{Introductory lectures on equivariant cohomology} Annals of Math. Studies 204, Princeton Univ. Press, 2020.
\end{thebibliography}
\end{document}